\documentclass[11 pt]{amsart}

\usepackage[margin=1.3in]{geometry}

\usepackage{biblatex}
\usepackage{xcolor}
\usepackage{amssymb}
\usepackage{amsmath}
\usepackage{amsthm}
\usepackage{amsfonts}
\usepackage{amsxtra}
\usepackage[unicode=true, pdfusetitle,
 bookmarks=true,bookmarksnumbered=false,
 breaklinks=false,
 backref=false,
 colorlinks=true,
 linkcolor=black,
 citecolor=teal,
 urlcolor=blue,
 final
]{hyperref}

\newcommand{\newrefformat}[2]{}

\theoremstyle{plain}   
\newtheorem{thm}{Theorem}

\newtheorem{cor}{Corollary}

\newtheorem{lem}{Lemma}[section]

\title{An improved upper bound for Tuza's conjecture via 2-colorable triangle families}

\author{Lixing Yi}
\begin{document}
\begin{abstract}
Tuza's conjecture states that for any graph $G$, the minimum size of a triangle transversal $\tau(G)$ is at most twice the maximum size of a set of edge-disjoint triangles $\nu(G)$. In this note, we prove $\tau(G) \leq \frac{63}{22}\nu(G)$, improving the previous bound $\tau(G) \leq \frac{66}{23}\nu(G)$ established by Haxell in 1999. The key observation is that for a ``2-colorable'' family of triangles $\mathcal{F}$, where each triangle has two blue edges and one red edge, we can obtain $\tau(\mathcal{F})\leq (1+\sqrt{3})\nu(\mathcal{F})$.
\end{abstract}
\maketitle

\section{Introduction}
A well-known conjecture of Tuza \cite{TUZA} asserts that in any simple graph $G$, the minimum size of a triangle transversal $\tau(G)$ and the maximum number of edge-disjoint triangles $\nu(G)$ satisfy the inequality $\tau(G) \leq 2\nu(G)$. While the conjecture has been proven for several classes of graphs, such as planar graphs \cite{planar}, graphs with bounded treewidth \cite{treewidth}, threshold graphs \cite{Thres}, etc., the general case remains open. 

The trivial upper bound is $\tau(G) \leq 3\nu(G)$. In 1999, Haxell \cite{HAXELL} improved the general bound to $\tau(G) \leq \frac{66}{23}\nu(G)$, which has stood as the best general bound. In the same paper, Haxell noted that the bound can be improved to $\tau(G)\leq\frac{1+\sqrt{481}}{8}\nu(G)$, but the proof was omitted.

In this note, we observe that a triangle family in Haxell's proof is ``2-colorable''. We prove the bound $\tau(\mathcal{F})\leq (1+\sqrt{3})\nu(\mathcal{F})$ for any 2-colorable triangle family $\mathcal{F}$. Combining Haxell's original results, we prove $\tau(G) \leq \frac{63}{22}\nu(G)$ for any $G$. Using the irrational constant $1+\sqrt{3}$ yields the slightly sharper bound $\tau(G) \leq \frac{162+4\sqrt{3}}{59}\nu(G)$.

\section{A bound on 2-colorable triangle families}

A family of triangles $\mathcal{F}$ in a graph $G$ is \textit{2-colorable} if we can color the edges of $G$ red and blue, such that each triangle in $\mathcal{F}$ has two blue edges and one red edge. For the rest of the section, ``triangle'' refers to a triangle in the family, unless otherwise stated. We call a family of triangles independent if they are edge-disjoint.

Given $\mathcal{F}$, let $\mathcal{B}$ be a maximum-size independent family of triangles in $\mathcal{F}$. Denote $\nu(\mathcal{F})=|\mathcal{B}|$. Let $\tau(\mathcal{F})$ be the minimum size of a triangle transversal of $\mathcal{F}$, i.e. an edge set intersecting every triangle of $\mathcal{F}$.

\begin{thm}\label{thm}
    For a 2-colorable triangle family $\mathcal{F}$, we have $\tau(\mathcal{F})\leq (1+\sqrt{3})\nu(\mathcal{F})$.
\end{thm}
To prove Theorem \ref{thm}, we need two lemmas and some definitions. 

We let $E(\mathcal{B})$ denote the edge set of $\mathcal{B}$. We say a triangle is type $(\mathcal{B},i)$ if it shares exactly $i$ edge(s) with $E(\mathcal{B})$. Let $\textcolor{blue}{E_B}$ be the set of blue edges in $E(\mathcal{B})$. Let $\textcolor{red}{E_R}$ be the set of red edges in $E(\mathcal{B})$, and let $\mathcal{F}_R$ be the set of type $(\mathcal{B},1)$ triangles that intersect $E(\mathcal{B})$ at $\textcolor{red}{E_R}$. Let $\mathcal{B}_1$ be a maximum independent subset of $\mathcal{F}_R$.

The following lemma constructs a transversal that exploits the local structure of $\mathcal{B}_1$. Its proof is almost identical to the proof of \cite[Lemma 1]{HAXELL}, and we include it for completeness.
\begin{lem}\label{lem1}
    $\tau(\mathcal{F})\leq 3\nu(\mathcal{F}) - |\mathcal{B}_1|.$
\end{lem}
\begin{proof}
Since $\mathcal{B}_1$ is independent, each triangle $T$ in $\mathcal{B}_1$ has a unique red edge and intersects a unique triangle $T'$ in $\mathcal{B}$. Let $\mathcal{B}'\subseteq \mathcal{B}$ denote the family $\{T':T\in\mathcal{B}_1\}$, then $|\mathcal{B}'|=|\mathcal{B}_1|$. Now, take an intersecting pair $T\in \mathcal{B}_1$ and $T'\in \mathcal{B}'$. The subgraph of $G$ formed by $E(T) \cup E(T')$ is a copy of $K_4$ minus an edge. Let $e(T)$ denote the red edge shared by $T$ and $T'$, and if the remaining edge of the $K_4$ also exists in $G$, let it be denoted by $e'(T)$. 

Take any triangle $P$ that is edge-disjoint from $\mathcal{B}\setminus\mathcal{B}'$, we will show it must share an edge with both $T$ and $T'$ for some $T\in\mathcal{B}_1$. Since $P$ must intersect $\mathcal{B}$ but is edge-disjoint from $\mathcal{B}\setminus\mathcal{B}'$, it must intersect $\mathcal{B}'$. Let $\mathcal{B}_P'\subseteq \mathcal{B}'$ be the subset of $\mathcal{B}'$ that intersects $P$, and let $\mathcal{B}_P$ be the corresponding subset in $\mathcal{B}_1$. Consider $(\mathcal{B}\setminus \mathcal{B}_P')\cup \mathcal{B}_P$, which is a maximum-size independent subfamily of $\mathcal{F}$. Hence $P$ must intersect $\mathcal{B}_P$, thus share an edge with some $T\in\mathcal{B}_P$. $P$ already shares an edge with $T'\in\mathcal{B}'_P$, so $P$ shares an edge with both $T$ and $T'$. Then $P$ must contain $e(T)$ or $e'(T)$.

Then $C = E[\mathcal{B}\setminus\mathcal{B}'] \cup \{e(T): T\in \mathcal{B}_1\} \cup \{e'(T): T\in \mathcal{B}_1 \text{ if $e'(T)$ exists}\}$ is a transversal of $\mathcal{F}$ and its size is at most $3\nu(\mathcal{F}) - |\mathcal{B}_1|$.
\end{proof}
The next lemma leverages the 2-coloring of $\mathcal{F}$ by using the blue edges to efficiently cover most triangles in $\mathcal{F}$.
\begin{lem}\label{lem2}
    $\tau(\mathcal{F})\leq |\textcolor{blue}{E_B}|+ \tau(\mathcal{F}_R)=2\nu(\mathcal{F})+ \tau(\mathcal{F}_R).$
\end{lem}
\begin{proof}
    We will construct a triangle transversal for $\mathcal{F}$. First, add $\textcolor{blue}{E_B}$, the blue edges of $E(\mathcal{B})$ to the transversal. Since all triangles intersect $E(\mathcal{B})$, the only triangles not covered by $\textcolor{blue}{E_B}$ must intersect $E(\mathcal{B})$ at $\textcolor{red}{E_R}$, and they must be type $(\mathcal{B},1)$ since each triangle has one red edge. Hence the triangles not covered are precisely $\mathcal{F}_R$. Let $C$ be a minimum transversal of $\mathcal{F}_R$, then $\textcolor{blue}{E_B}\cup C$ is a transversal of $\mathcal{F}$. $|C|=\tau(\mathcal{F}_R)$, and since each triangle in $\mathcal{B}$ has two blue edges, $|\textcolor{blue}{E_B}|=2\nu(\mathcal{F})$. Hence $\tau(\mathcal{F})\leq |\textcolor{blue}{E_B}\cup C|= 2\nu(\mathcal{F})+\tau(\mathcal{F}_R).$
\end{proof}

Now we can prove Theorem \ref{thm}.
\begin{proof}
    We will induct on $|\mathcal{B}|=\nu(\mathcal{F})$. When $\nu(\mathcal{F})=|\mathcal{B}|=0$, $\mathcal{F}$ is empty, so $\tau(\mathcal{F})\leq (1+\sqrt{3})\nu(\mathcal{F})$ vacuously holds. Now suppose Theorem \ref{thm} holds for all $\mathcal{F}'$ with $\nu(\mathcal{F}')<n$, and say $\nu(\mathcal{F})=n$. For the inductive step, either $|\mathcal{B}_1|=\nu(\mathcal{F})$, then  $\tau(\mathcal{F})\leq 3\nu(\mathcal{F}) - |\mathcal{B}_1|=2\nu(\mathcal{F})$ by Lemma \ref{lem1} and we are done, or $|\mathcal{B}_1|<\nu(\mathcal{F})$. Then $\nu(\mathcal{F}_R)=|\mathcal{B}_1|<n$, and since $\mathcal{F}_R$ is 2-colorable, we can apply the inductive hypothesis on $\mathcal{F}_R$ and get $\tau(\mathcal{F})\leq |\textcolor{blue}{E_B}|+\tau(\mathcal{F}_R)\leq 2\nu(\mathcal{F}) +(1+\sqrt{3})|\mathcal{B}_1|$ by Lemma \ref{lem2}. Combining this with Lemma \ref{lem1} yields \begin{align*}
(1+\sqrt{3})\tau(\mathcal{F})+\tau(\mathcal{F})&\leq (1+\sqrt{3})(3\nu(\mathcal{F}) - |\mathcal{B}_1|)+2\nu(\mathcal{F}) +(1+\sqrt{3})|\mathcal{B}_1|\\
(2+\sqrt{3})\tau(\mathcal{F})&\leq (5+3\sqrt{3})\nu(\mathcal{F})\\
\tau(\mathcal{F})&\leq (1+\sqrt{3})\nu(\mathcal{F})
    \end{align*}
\end{proof}
\section{An improved general bound on Tuza's conjecture}
Theorem \ref{thm} directly improves the best known general bound $\tau(G)\leq\frac{66}{23} \nu(G)$ and the sketched bound $\tau(G)\leq\frac{1+\sqrt{481}}{8}\nu(G)\approx 2.8665 \nu(G)$ given by Haxell \cite{HAXELL}.
\begin{cor}\label{cor}
    For a graph $G$, we have $\tau(G)\leq\frac{63}{22} \nu(G)\approx 2.8636 \nu(G)$.
\end{cor}
For completeness, we restate the necessary definitions and the four lemmas in \cite{HAXELL}. The reader should refer to \cite{HAXELL} for their proofs.

Given a graph $G$, let $\mathcal{B}$ be a maximum-size independent family of triangles. Let $\mathcal{B}_1$ be a maximum-size independent family of type $(\mathcal{B},1)$ triangles. Let $G'=G\setminus E(\mathcal{B}_1)$. Let $\mathcal{B}_2$ be a maximum-size independent family of type $(\mathcal{B},2)$ triangles in $G'$. Let $\mathcal{B}'$ be a maximum-size independent family of triangles in $G'$, subject to the condition that $|E(\mathcal{B}')\setminus E(\mathcal{B})|\geq |\mathcal{B}_2|$. Let $\mathcal{S}$ be the family of type $(\mathcal{B}',1)$ triangles that share exactly one edge with $E(\mathcal{B}')\setminus E(\mathcal{B})$, and $\mathcal{B}_1'$ a maximum-size independent subset of $\mathcal{S}$. We have: 

\begin{lem}\cite[Lemma 1]{HAXELL}
     $\tau(G)\leq 3\nu(G)-|\mathcal{B}_1|$.
\end{lem}
\begin{lem}\cite[Lemma 2]{HAXELL}
     $\tau(G)\leq \frac{3}{2}\nu(G)+\frac{5}{2}|\mathcal{B}_1|+2|\mathcal{B}_2|$.
\end{lem}
\begin{lem}\cite[Lemma 3]{HAXELL}\label{lem3}
     $\tau(G)\leq 3\nu(G)-|\mathcal{B}_1'|$.
\end{lem}
\begin{lem}\cite[Lemma 4]{HAXELL}\label{lem4}
     $\tau(G)\leq 3\nu(G)+3|\mathcal{B}_1'|-|\mathcal{B}_2|$.
\end{lem}
Now we can prove Corollary \ref{cor}.
\begin{proof}
    In Haxell's proof of Lemma \ref{lem4}, a triangle transversal of $G$ was constructed by first adding $E(\mathcal{B}_1)$ and $E(\mathcal{B}')\cap E(\mathcal{B})$. Note that $|E(\mathcal{B}_1)\cup (E(\mathcal{B}')\cap E(\mathcal{B}))|\leq 3\nu(G)-|\mathcal{B}_2|$. The triangles left to cover are precisely $\mathcal{S}$. Haxell used $E(\mathcal{B}_1')$ as a transversal of $\mathcal{S}$, which contributed the term $3|\mathcal{B}_1'|$ in Lemma \ref{lem4}. Notice that $\mathcal{S}$ is 2-colorable: we can color edges of $G$ by ``In $\mathcal{B}'$\ '' and ``Not in $\mathcal{B}'$\ '', then since all triangles in $\mathcal{S}$ are type $(\mathcal{B}',1)$, $\mathcal{S}$ has an edge 2-coloring, where each triangle has exactly one edge colored ``In $\mathcal{B}'$\ ''. Hence we can apply Theorem \ref{thm} and bound a minimum transversal of $\mathcal{S}$ by $(1+\sqrt{3})|\mathcal{B}_1'|$. Here, we will use the weaker but simpler rational bound $\tau(\mathcal{S})\leq\frac{11}{4}|\mathcal{B}_1'|$. Unioning a minimum $\mathcal{S}$ transversal with $E(\mathcal{B}_1)\cup (E(\mathcal{B}')\cap E(\mathcal{B}))$ yields a transversal of $G$, hence $\tau(G)\leq 3\nu(G)+\tau(\mathcal{S})-|\mathcal{B}_2|\leq 3\nu(G)+\frac{11}{4}|\mathcal{B}_1'|-|\mathcal{B}_2|.$

    With the improved version of Lemma \ref{lem4}, we combine the rest of the lemmas:
    \begin{align*}
        \frac{5}{2}\tau(G)+\tau(G)+\frac{11}{2}\tau(G)+2\tau(G)&\leq \frac{5}{2}(3\nu(G)-|\mathcal{B}_1|)+\frac{3}{2}\nu(G)+\frac{5}{2}|\mathcal{B}_1|+2|\mathcal{B}_2|\\
        &+ \frac{11}{2}(3\nu(G)-|\mathcal{B}_1'|) + 2(3\nu(G)+\frac{11}{4}|\mathcal{B}_1'|-|\mathcal{B}_2|)\\
        \tau(G)&\leq \frac{63}{22}\nu(G)
    \end{align*}
\end{proof}
In fact, directly applying the constant $1+\sqrt{3}$ in Theorem \ref{thm} yields the sharper bound $\tau(G)\leq\frac{162+4\sqrt{3}}{59} \nu(G)$, at the cost of elegance. An obvious next step, therefore, is to tighten the bound for 2-colorable triangle families in Theorem \ref{thm}, since any improvement will immediately improve the general bound for Tuza's conjecture. However, note that the triangles in $K_4$ are 2-colorable, and $\tau(K_4)=2\nu(K_4)$, so Tuza's conjectured constant 2 is best possible in Theorem \ref{thm}. Thus improving the bound for 2-colorable triangle families can at best improve the general constant to 54/19.
\section*{Acknowledgements}
The author thanks Fan Chung for her encouragement and helpful discussions.

\section*{Disclosure of AI use}
The results were developed by the author independently without AI tools. Generative LLMs (ChatGPT and Claude) were used for reviewing and editing the manuscript only. The author takes full responsibility for the content of the article.
\printbibliography
\end{document}